\documentclass[11pt]{amsart}
\usepackage{amsfonts,amssymb,amsthm,amsmath,amsxtra,amscd,verbatim,eucal}
\usepackage[all]{xy}
\usepackage[dvips]{graphics}

\theoremstyle{plain}
\newtheorem{thm}{Theorem}[section]
\newtheorem{lem}[thm]{Lemma}
\newtheorem{prop}[thm]{Proposition}
\newtheorem{cor}[thm]{Corollary}

\theoremstyle{definition}
\newtheorem{defi}[thm]{Definition}

\newtheorem{eg}[thm]{Example}

\theoremstyle{remark}
\newtheorem{rmk}[thm]{Remark}

\def\Z{{\mathbb Z}}
\def\N{{\mathbb N}}

\def\C{{\mathbb C}}
\def\A{{\mathbb A}}
\def\R{{\mathbb R}}
\def\Q{{\mathbb Q}}
\def\P{{\mathbb P}}

\def\O{\mathcal{O}}

\def\J{\mathcal{J}}

\def\aa{\mathfrak{a}}

\def\mm{\mathfrak{m}}

\def\a{\alpha}
\def\b{\beta}
\def\g{\gamma}
\def\d{\delta}

\def\e{\eta}

\def\n{\nu}

\def\p{\pi}

\def\s{\sigma}
\def\t{\tau}

\def\D{\Delta}

\def\S{\Sigma}
\def\Om{\Omega}

\def\.{\cdot}
\def\^{\widehat}
\def\~{\widetilde}
\def\o{\circ}
\def\ov{\overline}

\def\surj{\twoheadrightarrow}
\def\inj{\hookrightarrow}
\def\de{\partial}

\renewcommand{\and}{ \quad \text{ and } \quad }
\renewcommand{\for}{ \quad \text{ for } \ }
\newcommand{\fall}{ \quad \text{ for all } \ }
\newcommand{\where}{ \quad \text{where} \ }

\newcommand{\lra}{\longrightarrow}

\newcommand{\cont}{\operatorname{Cont}}
\let \cedilla =\c

\newcommand{\sTm}{{\Spec K[t]/(t^{m+1})}}
\newcommand{\sT}{{\Spec K[[t]]}}
\newcommand{\sing}{{\operatorname{Sing}}}

\DeclareMathOperator{\codim} {codim}

\DeclareMathOperator{\im}  {im}

\DeclareMathOperator{\Gr} {Gr}

\DeclareMathOperator{\Proj} {Proj}
\DeclareMathOperator{\Spec} {Spec}

\DeclareMathOperator{\Sing} {Sing}

\DeclareMathOperator{\Bl} {Bl}

\DeclareMathOperator{\Ex} {Ex}

\DeclareMathOperator{\ord} {ord}
\DeclareMathOperator{\reg} {reg}

\DeclareMathOperator{\val} {val}

\DeclareMathOperator{\Cont} {Cont}

\DeclareMathOperator{\orb} {orb}

\begin{document}

\title{Divisorial valuations via arcs}


\author{Tommaso de Fernex}
\address{Department of Mathematics, University of Utah, 155 South 1400
East, Salt Lake City, UT 48112-0090, USA}
\email{{\tt defernex@math.utah.edu}}

\author{Lawrence Ein}
\address{Department of Mathematics, University of Illinois at Chicago,
851 S. Morgan St., M/C. 249, Chicago, IL 60607-7045, USA}
\email{{\tt ein@math.uic.edu}}

\author{Shihoko Ishii}
\address{Department of Mathematics, Tokyo Institute of Technology,
Oh-Okayama, Meguro, 152-8551, Tokyo, Japan}
\email{{\tt shihoko@math.titech.ac.jp}}

\thanks{The research of the first autho was partially supported by
NSF grants DMS 0111298 and DMS 0548325,
and by MIUR National Research Project
``Geometry on Algebraic Varieties" (Cofin 2004).
The research of the second author was partially supported
by NSF grant DMS 0200278. The research of the third author
was partially supported by Grant-In-Aid of Ministry of
Science and Education in Japan.}


\maketitle

\begin{abstract}
This paper shows a finiteness property of
a divisorial valuation in terms of arcs.
First we show that every divisorial valuation
over an algebraic variety corresponds to an
irreducible closed subset of the arc space.
Then we define the codimension for this subset
and give a formula of the codimension in terms
of  ``relative Mather canonical class''.
By using this subset, we prove that a divisorial
valuation is determined by assigning the values of
finite functions. We also have a criterion for a divisorial
valuation to be a monomial valuation by assigning
the values of finite functions.
\end{abstract}

\section*{Introduction}

Let $X$ be a complex algebraic variety of dimension $n \ge 1$.
An important class of valuations of the function field $\C(X)$
of $X$ consists of {\it divisorial valuations}. These are
valuations of the form
$$
v  =  q\val_E  \colon \C(X)^* \lra \Z,
$$
where $E$ is a prime divisor on a normal variety $Y$ equipped with a
proper, birational morphism $f \colon Y \to X$,
$q = q(v)$ is a positive integer number
called the {\it multiplicity} of $v$, and for every $h \in  \C(X)^*$
that is regular at the generic point of $f(E)$,
$\val_E(h) := \ord_E(h\o f)$ is the order of vanishing of $h\o f$ at the
generic point of $E$.
It is a theorem of Zariski that the set of valuation rings associated
with this class of valuations coincides with the set of
discrete valuation rings $(R,\mm_R)$ of $\C(X)$
with ${\rm trdeg}(R/\mm_R : \C) = n-1$. Thanks to Hironaka's
resolution of singularities, divisorial valuations had
acquired a fundamental role in singularity theory.

Since the publication of the influential paper of Nash \cite{Nash}
and the introduction of motivic integration
(see, e.g., \cite{Kon,DL,Bat}), it become apparent
the close link between certain invariants of singularities related to
divisorial valuations and the geometry of arc spaces.
This link was first explored by Musta\c{t}\v{a}
in \cite{must01,must02}, and then further
studied in \cite{ELM,i2,i3}. In particular,
when the ambient variety $X$ is smooth,
it is shown in \cite{ELM} how one can
reinterpret invariants such as multiplier ideals
in terms of properties of certain subsets
in the space of arcs $X_\infty$ of $X$.

The main purpose of this paper is to extend the results
of \cite{ELM} to arbitrary varieties and to employ
such results towards the characterization of divisorial
valuations by evaluation against finite numbers of functions.

Naturally associated to the valuation $v = q\val_E$,
there is a subset
$$
W(v) = W(E,q) \subset X_\infty,
$$
constructed as follows.
We can assume without loss of generality that both $Y$ and $E$ are smooth.
Then $W(v)$ is defined as the closure of the image,
via the natural map $Y_\infty \to X_\infty$,
of the set of arcs on $Y$ with order of contact along $E$
equal to $q$. It turns out that
$$
v = \val_{W(v)},
$$
where for every irreducible closed subset $C\subset X_{\infty}$
that is not contained in the arc space of any
proper closed subvariety of $X$, we define a valuation
$\val_C\colon \C(X)^* \to \Z$ by taking the order of vanishing
(or polarity) along the generic point of $C$
(see Definition~\ref{defi:val-along-fat-sets}).
Although the subsets of $X_\infty$ of the form $W(v)$
are quite special, there is a much larger class
of subsets $C \subseteq X_\infty$, such that
the associated valuation $\val_C$ is a divisorial valuation.
These sets are called {\it divisorial sets}.
It was shown in \cite{i3} that sets of the form
$W(v)$ are maximal (with respect to inclusion)
among all divisorial sets defining the same valuation;
for this reason they are called {\it maximal  divisorial sets}.

Other important classes of subsets of $X_\infty$
are those consisting, respectively, of {\it (quasi)-cylinders}
and {\it contact loci}. The valuations
associated to irreducible components of such sets are called
{\it cylinder valuations} and {\it contact valuations},
respectively.
As it turns out, maximal divisorial sets belong to these
classes of sets. Generalizing to singular varieties
some results of \cite{ELM}, we obtain the following properties:
\begin{enumerate}
\item
every divisorial valuation is a
cylinder valuation (Theorem~\ref{thm:valuations});
\item
every  cylinder  valuation is a contact valuation (Proposition~\ref{cyl-con});
\item
every  contact valuation is a
divisorial valuation (Proposition~\ref{component}).
\end{enumerate}
Geometrically, the correspondence between divisorial valuations
and cylinder valuations is constructed by associating
to any divisorial valuation
$v$ the subset $W(v) \subset X_\infty$.
The fact that $W(v)$ is a quasi-cylinder shows that $v$, being
equal to $\val_{W(v)}$, is a cylinder valuation.
The other direction of the correspondence
is more elaborate: starting from a cylinder valuation $\val_C$,
one first realizes $\val_C$ as a contact valuation,
that is, a valuation determined by an irreducible component
of a contact locus, the definition of which involves certain
conditions on the order of contact along some subscheme of $X$.
The divisor $E$ is then extracted by a suitable
weighted blowup on a log resolution of this subscheme,
and $q$ is determined by the numerics involved in the construction.

A key point in the proof of these assertions,
as well as a fundamental property for many applications,
is a codimension formula for the maximal
divisorial set $W(v)$ associated to a divisorial
valuation $v$. When $X$ is a smooth variety,
it was shown in \cite{ELM} that
$$
\codim\big(W(v),X_\infty\big) = k_v(X) + q(v),
$$
where $k_v(X) := v(K_{Y/X})$ is the
discrepancy of $X$ along $v$ and $q(v)$ is the multiplicity of $v$.
The definition of discrepancy needs to be modified when $X$ is singular.
Given an arbitrary variety $X$, we take a resolution of singularities
$f \colon Y \to X$ that factors through the Nash blowup
$$
\n \colon \^X \lra X,
$$
and define the {\it relative Mather canonical divisor} $\^K_{Y/X}$
of $f$ (see Definition~\ref{defi:Nash-blowup-Mather-canonical-class}).
This divisor, which is defined in total generality, is always an
effective integral divisor, and it coincides with
$K_{Y/X}$ when $X$ is smooth (the two divisors
are in general different for $\Q$-Gorenstein varieties).
The relative Mather canonical divisor plays a fundamental
role in the geometry of arc spaces
and the change-of-variables formula in motivic integration.
Defining $\^k_v(X) := v(\^K_{Y/X})$,
the codimension formula of \cite{ELM}
generalizes to arbitrary varieties as follows.

\begin{thm}\label{thm:intro:valuations}
With the above notation, we have
$$
\codim\big(W(v),X_\infty\big) =  \^k_v(X) + q(v).
$$
\end{thm}

Using the interpretation of divisorial valuations via arc spaces,
we show that each divisorial valuation $v$
is characterized by its values $v(f_i)=v_i$ on a
finite number of functions $f_i$.
More precisely, we have the following result.

\begin{thm}\label{thm:intro:finite}
Suppose that $X = \Spec A$ is an affine variety,
and let $v$ be a divisorial valuation over $X$.
Then there exists elements $f_1,\dots,f_r \in A$ and
$v_1,\dots,v_r \in \N$ such that for every $f \in A\setminus\{0\}$
$$
v(f) =
\min\{ v'(f) \mid \text{$v'$ is a divisorial valuation
such that $v'(f_i) = v_i$} \}.
$$
\end{thm}

Theorem~\ref{thm:intro:finite} is obtained by determining
functions $f_i$ and numbers $v_i$ such that
\begin{equation}\label{eq:inclusion}
W(v)\subset \overline{\bigcap_{i=1}^r\cont ^{v_{i}}(f_{i})  },
\end{equation}
with equality holding off a set contained in $(\Sing X)_{\infty}$.
This is the content of Theorem~\ref{finite}.
A similar result can be obtained using MacLane's results
from \cite{ML1} (see Remark~\ref{rmk:MacLane});
it would be interesting to further investigate the
connection between MacLane's key polynomials
and the functions $f_i$ determined in the
proof of the above results, and to study properties that the
first ones satisfy with respect to the geometry of arc spaces.

In the case of monomial valuations
on toric varieties, we apply Theorem~\ref{thm:intro:valuations}
to give a precise characterization in terms of a system of parameters.
For simplicity, we present here the result in the
special case when the toric variety is equal to $\C^n$.

\begin{thm}\label{thm-intro:monomial-smooth}
Let $v$ be a divisorial valuation of
$\C(x_1,\dots,x_n)$, centered at the origin $0$ of $X = \C^n$.
Assume that there are positive integers $a_1,\dots,a_n$ such that
\begin{equation}\label{eq-intro:monomial-smooth}
v(x_i)  \ge  a_i \and
\sum a_i \ge k_v(X) + q(v).
\end{equation}
Then $v$ is a monomial valuation determined by the weights $a_i$
assigned along the parameters $x_i$, $q(v) = \gcd(a_1,\dots,a_n)$,
and equalities hold in both formulas in \eqref{eq-intro:monomial-smooth}.
\end{thm}

The more general version of this result, holding for
arbitrary singular toric varieties,
requires some additional notation, and is given
in Theorem \ref{thm-intro:monomial-singular}.

The authors express their hearty thanks to Bernard Teissier
for his  helpful advice, and
would like to thank Mircea Musta\c{t}\v{a} for many helpful
discussions, in particular about the proof of
Theorem~\ref{thm-intro:monomial-singular}.

\section{The relative Mather canonical divisor}\label{section:mather-canonical-divisor}

In this preliminary section, we review the construction
of the Nash blowup and define a generalization
of the relative canonical divisor to certain resolutions
of arbitrarily singular varieties. This will
give a geometric interpretation of a certain ideal
sheaf that governs the dimension of fibers of
maps between arc space, and consequently the change-of-variable formula
in motivic integration.

We start with an arbitrary complex variety $X$ of dimension $n$.
Note that the projection
$$
\p \colon \P_X(\wedge^n \Om_X) \lra X
$$
is an isomorphism over the smooth locus $X_{\reg} \subseteq X$.
In particular, we have a section $\s \colon X_{\reg} \to \P_X(\wedge^n \Om_X)$.

\begin{defi}\label{defi:Nash-blowup-Mather-canonical-class}
The closure of the image of the section $\s$ is the {\it Nash blowup} of $X$,
and is denoted by $\^X$:
$$
\xymatrix@C=5pt{
& \P_X(\wedge^n \Om_X) \ar[d]^\p
& \supseteq \; \^X := \ov{\s(X_{\reg})} \quad\quad\quad  \\
X_{\reg} \ar@{^{(}->}[r] \ar@/^15pt/[ur]^(.32)\s & X.
}
$$
The Nash blowup $\^X$ comes equipped with the morphism
$$
\n := \p|_{\^X} \colon \^X \lra X
$$
and the line bundle
$$
\^K_X := \O_{\P_X(\wedge^n \Om_X)}(1)|_{\^X}.
$$
We call this bundle the {\it Mather canonical line bundle} of $X$.
\end{defi}

\begin{rmk}
If $X$ is smooth, then $\^X = X$ and $\^K_X$ is just the
canonical line bundle of $X$.
\end{rmk}

The original definition of Nash blowup is slightly different.
Assuming the existence of an embedding $X \inj M$ into a manifold $M$,
one can consider the Grassmann bundle $G(\Om_M,n)$ over $M$ of rank $n$
locally free quotients sheaves of $\Om_M$. The map
$$
x \mapsto \big((\Om_X)_x \twoheadleftarrow (\Om_M)_x \big) \in G(\Om_M,n)_x,
$$
defined for every smooth point $x$ of $X$, gives
a section $\~\s \colon X_{\reg} \to G(\Om_M,n)$.
Then one takes the closure of the image of this section in $G(\Om_M,n)$.
As we prove next, the two constructions agree.

\begin{prop}
Keeping the above notation, let $\~X$ denote the closure of
$\~\s(X_{\reg})$ in $G(\Om_M,n)$, and let
$\~\n \colon \~X \to X$ the induced morphism.
Then $\~X = \^X$ and $\~\n = \n$.
\end{prop}

\begin{proof}
We have exact sequences
$$
\xymatrix{
0 \ar[r] & S|_{\~X} \ar[r]^i & \Om_M|_{\~X} \ar[r]^p \ar@{=}[d]
& Q|_{\~X} \ar[r] & 0 \\
0 \ar[r] & \ker(q) \ar[r]^j & \Om_M|_{\~X} \ar[r]^q & \~\n^*\Om_X \ar[r] & 0,
}
$$
where $ \Om_M|_{\~X}$ is $\varpi^*\Om_M|_{\~X}$ for the projection
$\varpi\colon G(\Om_M, n)\to M$ and it also coincides with $\~\n^*(\Om_M|_X)$.
The top row is the restriction of the universal sequence of the
Grassmann bundle $G(\Om_M,n)$, and the second is the pull back of the
sequence of differentials determined by the inclusion of $X$ in $M$.
Over the smooth locus of $X$ we have
$$
S|_{\~\s(X_{\reg})} = \ker(q)|_{\~\s(X_{\reg})}.
$$
Then, since $\ker(q)$ is torsion free and the
top sequence has a local splitting, the inclusion $j$ factors
through $i$ and an inclusion $\ker(q) \inj S|_{\~X}$.
 Therefore
$p$ factors through $q$ and a surjection $\~\n^*\Om_X \surj Q|_{\~X}$.
After taking wedges, we obtain a commutative diagram of surjections
\begin{equation}\label{eq:facoring-surjections}
\xymatrix@C=20pt@R=20pt{
\wedge^n \Om_M|_{\~X} \ar[dr] \ar[r] & \wedge^n Q|_{\~X} \\
&\wedge^n \~\n^*\Om_X \ar[u]
}.
\end{equation}
Now we consider the inclusion over $\~X$
$$
G(\Om_M|_{\~X},n) \inj \P_{\~X}(\wedge^n \Om_M|_{\~X})
$$
given by Pl\"ucker embedding.
The factorization \eqref{eq:facoring-surjections}
implies that the image of $\~X$ under this embedding
is contained in $\P_{\~X}(\wedge^n \~\n^* \Om_X)$,
when the latter is viewed as a subvariety of
$\P_{\~X}(\wedge^n \Om_M|_{\~X})$ via the natural embedding.
Therefore, by the compatibility of Pl\"ucker embeddings,
  the image of $\~X$ under the embedding
$$
G(\Om_M|_{X},n) \inj \P_{X}(\wedge^n \Om_M|_{X})
$$
is contained in $\P_X(\wedge^n\Om_X)$.
Then, restricting over the regular locus of $X$, we have
$$
\~X \cap \p^{-1}(X_{\reg}) = \p^{-1}(X_{\reg}) = \^X \cap \p^{-1}(X_{\reg}).
$$
Since $\~X$ and $\^X$ are irreducible varieties, both
surjecting onto $X$, we conclude
that $\~X = \^X$. The equality $\~\n = \n$ follows by the fact that the
construction is carried over $X$.
\end{proof}

\begin{rmk}
The original construction of the Nash blowup using Grassmann bundles
can be given without using (or assuming) any embedding,
by considering the Grassmann bundle $G(\Om_X,n)$ on $X$
of rank $n$ locally free quotient sheaves of $\Om_X$.
Notice that $G(\Om_X,n)$ and $\P_X(\wedge^n \Om_X)$ agree
over the smooth locus of $X$.
\end{rmk}

\begin{rmk}
Mather used the above construction to propose
a generalization to singular varieties
of the notion of Chern classes of manifolds, by considering the class
$\~\n_* \big(c(Q^*|_{\~X}) \cap [\~X]\big)$ in $A_*(X)$.
This is known as the {\it Mather-Chern class} of $X$.
The push-forward $\n_* \big(c_1(\^K_X^{-1}) \cap [\^X]\big)$ is equal to the
first Mather-Chern class of $X$.
\end{rmk}

Now consider any resolution of singularities $f \colon Y \to X$
factoring through the Nash blowup of $X$, so that we have a commutative
diagram
$$
\xymatrix{
Y \ar[r]_{\^f} \ar@/^15pt/[rr]^f & \^X \ar[r]_{\n} & X
}.
$$
\begin{defi}
  Let $\^K_{Y/X}$ be the divisor  supported on the exceptional divisor on $Y$
  and linearly equivalent to $K_Y- \^f^*\^K_X$.
  We call it the {\it relative Mather canonical divisor}
of $f$.
\end{defi}

\begin{prop}
The relative Mather canonical divisor $\^K_{Y/X}$ is an effective divisor and satisfies the relation:
$$df(f^*\wedge^n\Om_X)=\O_Y(-\^K_{Y/X})\cdot\wedge^n\Om_Y,$$
where $df: f^*\wedge^n \Om_X \to \wedge^n \Om_Y $ is the canonical homomorphism.
\end{prop}

\begin{proof}
By generic smoothness of $X$, the kernel of
the morphism $\n^*\wedge^n\Om_X \to \^K_X$ is torsion, hence,
pulling back to $Y$, we obtain a commutative diagram
$$
\xymatrix@C=20pt@R=20pt{
f^*\wedge^n \Om_X \ar[r]^{df} \ar[d] & \wedge^n \Om_Y \\
\^f^*\^K_X \ar[ur]_\d
}.
$$
Then, we have an effective divisor $D$ with the support on the exceptional divisor
such that
$$
\im(\d) = \O_Y(-D)\.\wedge^n\Om_Y.
$$
It follows that   $D$ is linearly equivalent to $K_Y- \^f^*\^K_X$,
therefore we obtain  $D=\^K_{Y/X}$.

 For the second statement, we should note that  $\n^*\wedge^n\Om_X \to \^K_X$
 is surjective as $\^K_X$ is relatively very ample with respect to $\p:\P_X(\wedge^n
 \Om_X)\to X$.
 This gives the surjectivity of  $f^*\wedge^n \Om_X \to \^f^*\^K_X $ and the second statement.
\end{proof}

Note that $\^K_{Y/X}$ is always an effective integral Cartier
divisor, and in particular it is in general different from
the relative canonical divisor defined in the $\Q$-Gorenstein case.
In fact, the following property holds.

\begin{prop}
Let $X$ be a normal and locally complete intersection variety,
and let $f \colon Y \to X$ be a resolution of singularities
factoring through the Nash blowup of $X$.
Then $\^K_{Y/X} = K_{Y/X}$ if and only if $X$ is smooth.
\end{prop}

\begin{proof}
It follows by \cite{EM1} that the difference $\^K_{Y/X} - K_{Y/X}$
is given by the vanishing order of the Jacobian ideal sheaf of $X$.
\end{proof}

\begin{defi}
For every prime divisor $E$ on $Y$, we define
$$
\^k_E(X) := \ord_E(\^K_{Y/X}),
$$
and call it the {\it Mather discrepancy} of $X$ along $E$.
More generally, if $v$ is a divisorial valuation over $X$,
then we can assume without loss of generality that
$v = q\val_E$ for a prime divisor $E$ on $Y$
and a positive integer $q$, and define
the {\it Mather discrepancy} of $X$ along $v$ to be
$$
\^k_v(X) := q\.\^k_E(X).
$$
If $X$ is smooth, then we denote
$k_E(X) := \^k_E(X)$ and $k_v(X) := \^k_v(X)$.
\end{defi}


\section{Contact loci in arc spaces and valuations}

\noindent

In this section we set up basic statements for
contact loci and divisorial valuations.

\begin{defi}
  Let  $ X $ be a scheme of finite type over $ \C$
and $K\supset \C$ a field extension.
  A morphism  $\alpha\colon \Spec K[[t]]\to X $
is called an {\it arc} of $ X $.
  We denote the closed point of $ \Spec K[[t]] $ by $ 0 $ and
  the generic point by $ \eta $.
  For $ m\in \N $, a morphism
$ \beta\colon \sTm \to X$ is called an {\it $ m $-jet} of $
  X $.
  Denote the space of arcs of $ X $ by $ X_{\infty} $ and the
  space of $ m $-jets of $ X $ by $ X_{m} $.
See \cite{must02} for more details.
\end{defi}

  The concept ``thin'' in the following is first introduced in
  \cite{ELM} and a ``fat arc'' is introduced and studied in \cite{i2}.

\begin{defi}
   Let $ X $ be a variety over $ \C$.
   We say that an arc $ \alpha\colon \sT  \to  X $ is {\it thin}
   if $ \alpha  $ factors through a proper closed subset of $ X $.
  An arc which is not thin is called a {\it fat arc}.
  An irreducible closed subset
$ C $ in $ X_{\infty} $  is called a {\it
  thin set} if  the generic point of $ C $ is thin.
  An irreducible closed subset  in $ X_{\infty} $
which is not thin  is called a {\it
  fat set}.

\end{defi}

 \begin{defi}\label{defi:val-along-fat-sets}
  Let $ \alpha\colon\sT\to X $ be a fat arc of a variety $ X $ and
  $ \alpha^*\colon \O_{X, \alpha(0)}\to K[[t]] $ be the  local homomorphism
  induced from $ \alpha $. Suppose that the induced morphism
  $\Spec K \to X$ is not dominant.
  By Proposition 2.5, (i) in  \cite{i2},
 $ \alpha^* $ is extended to the injective
  homomorphism of fields $ \alpha^*\colon\C(X)\to K((t)) $,
  where $ \C(X) $ is the rational function field of $ X $.
  Define a map $\val_{\alpha}\colon \C(X)\setminus\{0\}\to \Z $ by
  $$ \val_{\alpha}(f)=\ord_{t} (\alpha^*(f)). $$
  The function $\val_{\alpha} $ is a discrete valuation of $ \C(X) $.
  We call it the {\it  valuation corresponding to } $ \alpha $.
  If $\alpha$ is the generic point of an (irreducible and fat)
closed set
  $C$, the valuation $\val_{\alpha}$ is also denoted by
  $\val_C$, and is called the {\it valuation corresponding to } $C$.
  From now on, we denote $\ord_t\a^*(f)$ by $\ord_\a(f)$.
   A fat arc $ \alpha$ of $ X$ is called a {\it divisorial arc} if $
  \val_{\alpha} $ is a divisorial valuation over $ X $.
  A fat set $C$ is called a {\it divisorial set} if the valuation $\val_C$ is
  a divisorial valuation over $X$.
\end{defi}

\begin{rmk}
For every irreducible, fat set $C \subset X_\infty$
and every regular function $f$ on $X$, we have
  $  \val_C(f)=\min \{\ord_{\gamma}(f)\mid \gamma \in C\}$.
\end{rmk}

\begin{defi}
 Let \( \psi_{m}: X_{\infty}\to X_{m} \) be the canonical projection
 to the space of \( m \)-jets \( X_{m} \).
 A subset \( C\subset X_{\infty} \) is called a {\it cylinder} if
 there is a constructible set \( \Sigma\subset X_{m} \) for some \( m\in
 \N \) such that
 \[ C=\psi_{m}^{-1}(\Sigma). \]
\end{defi}

\begin{defi}[\cite{ELM}]
  For an ideal sheaf \( \aa   \) on a variety \( X \),
  we define
   \[  \cont^m(\aa  )= \{\alpha \in X_{\infty} \mid \ord_\a(\aa  )=
   m \} \]
      and
      \[  \cont^{\geq m}(\aa  )= \{\alpha \in X_{\infty} \mid \ord_\a(\aa  )
   \geq m \}. \]
     These subsets are called {\it contact loci} of the ideal \( \aa   \).
    The subset \( \cont^{\geq m}(\aa  ) \) is closed and
   \( \cont^m(\aa  ) \) is locally closed;
both are cylinders.
   If \( Z=Z(\aa  ) \), then \(  \cont^m(\aa  ) \) and
    \(  \cont^{\geq m}(\aa  ) \) are sometimes denoted by
   \(   \cont^{ m}(Z) \) and \(  \cont^{\geq m}(Z) \), respectively.
\end{defi}

\begin{defi}[\cite{ELM}]
  For a  simple normal crossing divisor \( E=\bigcup_{i=1}^sE_{i} \) on a
  non-singular variety \( X \), we introduce the {\it multi-contact locus}
  for a multi-index \( \nu=(\nu_1,\ldots \nu_s)\in \Z_{\geq 0}^s \):

   \[ \cont^{\nu}(E)=\{\alpha\in X_{\infty} \mid \ord_\a(I_{E_{i}})=
\nu_{i}\ \mbox{for\ every\ }
   i\}, \]
   where \(  I_{E_{i}} \) is the defining ideal of \( E_{i} \).
  If all intersections among the $E_i$ are irreducible,
  then each of these multi-contact loci \( \cont^{\nu}(E) \)
  is irreducible whenever it is   not empty.
\end{defi}

\begin{defi}
Let $f\colon Y\to X$ be a resolution of the singularities of $X$, and
suppose that $E$ is an irreducible smooth divisor on $Y$.
For any $q \in \Z_+$, we define
$$
W(E,q)=\overline{f_{\infty}(\cont^q(E))}
$$
  and call it a {\it maximal divisorial set}.
  For $v=q\val_E$ we denote sometimes $W(E,q) $ by $W(v)$.
  When we should clarify the space $X$ with $W(E,q)=W(v)\subset X_\infty$
we denote $W(E,q)$ by $W_X(E,q)$ or $W_X(v)$.
\end{defi}

\begin{rmk}
It follows by \cite[Proposition~3.4]{i3} that the above definition of
maximal divisorial set agrees with the one given in \cite[Definition 2.8]{i3}.
\end{rmk}

\begin{rmk}\label{rmk:I3}
The set $W(E,q)$ only depends on the valuation $v = q\val_E$,
and not on the particular model $Y$ we have chosen;
this justify the notation $W(v)$.
Moreover, let $g \colon X' \to X$ be a proper birational morphism
of normal varieties, and let $U \subset X'$ be an open subset intersecting
the center of $v$ on $X'$. We consider $v$ also as
a divisorial valuation over $X'$ and over $U$.
Since we can assume that $v = q\val_E$ for some smooth divisor $E$
on a resolution of $X'$,
it follows immediately from our definition of maximal divisorial set that
$W_X(v) = \overline{g_{\infty}(W_{X'}(v))}
= \overline{g_{\infty}(W_{U}(v))}$
(cf. \cite[Proposition~2.9]{i3}).
\end{rmk}

\begin{rmk}\label{rmk:maximality}
 The set $W(E,q)$ is a divisorial set corresponding to the
valuation $q\val_E$ and has
  the following ``maximality'' property: any divisorial
set $C$ with $\val_C=q\val_E$ is contained in
  $W(E,q)$ (see \cite{i3}).
\end{rmk}

The following is a generalization of a result of \cite{ELM}.

\begin{prop}
\label{component}
 Let \( X = \Spec A\) be an affine variety, and let \( \aa\subset A \)
be  a non-zero ideal.
 Then, for any \( m\in \N \),
 every fat irreducible component of  \( \cont^{\geq
 m}(\aa) \) is a maximal divisorial set.
\end{prop}

\begin{proof}
 Let \( \varphi\colon Y\to X \) be a log-resolution of \( \aa \), and write
 \( \aa\.\O_{Y}=\O_{Y}(-\sum_{i=1}^s r_{i}E_{i}) \), where
\( E=\bigcup_{i=1}^sE_{i} \) is a simple normal crossing divisor on $Y$.
By \cite[Theorem 2.1]{ELM}, we have
  \[ \cont^{\geq m}(\aa)\supset \bigsqcup_{\sum r_{i}\nu_{i}\geq m}
  \varphi_{\infty}(\cont^{\nu}(E)), \]
  where the complement in $ \cont^{\geq m}(\aa)$ of the above union is thin.

We claim that there are only finite number of maximal
divisorial sets
$ \overline{\varphi_{\infty}(\cont^{\nu}(E))}$'s
in $X_\infty$, for all possible values of $\nu$. This follows
by the following two facts:
  \begin{enumerate}
  \item[(i)]
We have  \( \nu\leq \nu' \) if and only if
\( \overline{\cont^{\nu}(E)}\supset \cont^{\nu'}(E) \),
where the partial order   \( \leq \) in \( \Z_{\geq 0}^s \) is defined by
  \[ (\nu_{1},\ldots,\nu_{s})\leq (\nu'_{1},\ldots,\nu'_{s}) \ \
  \mbox{if}\ \ \nu_{i}\leq \nu'_{i} \mbox{\ \ for\ all\ } i. \]

\item[(ii)]
The number of minimal elements of \( \{\nu\in \Z_{\geq 0}^s \mid
 \sum r_{i}\nu_{i}\geq m \}\) according to this order $\leq$ is finite.
   \end{enumerate}

  Then the maximal   \( \overline{\varphi_{\infty}(\cont^{\nu}(E)}) \)'s are
  the fat components of \( \cont^{\geq m}(\aa) \).
  By  \cite[Corollary 2.6]{ELM}, \( \overline{\cont^{\nu}(E)} \)'s are
  divisorial sets.
  Therefore, a fat irreducible component of \( \cont^{\geq m}(\aa) \)
  is a divisorial set.
 To show the maximality, let \( C \) be
a fat component of \( \cont^{\geq m}(\aa) \)
 and \( \alpha\in C \) the generic point.
 Let $\val_{\alpha}=q\val_F$.
 Then, it is clear that \( C \subset W(F,q) \) by the
 maximality of \(  W(F,q) \).
 For the opposite inclusion, take the generic point \( \beta\in
  W(F,q) \).
 Then it follows that \( \val_{\beta}=\val_{\alpha} \), which means that
 \( \ord_\b(f)=\ord_\a(f) \) for every \( f\in K(X) \).
 This gives \( \beta\in \cont^{\geq m}(\aa) \), and therefore
 \( W(F,q) \) is contained in a fat irreducible component
 of \(  \cont^{\geq m}(\aa)  \).
 In conclusion, \( C= W(F,q)\).
\end{proof}



\section{Codimension of a maximal divisorial set}

In this section we give an extension of the
formula on the codimension of maximal divisorial sets
established in \cite{ELM} to singular varieties.
Let $X$ be an arbitrary complex variety, and let $n = \dim X$.
Let $\J_X \subset \O_X$ be the Jacobian ideal sheaf of $X$.
In a local affine chart this ideal is defined as follows.
Restrict $X$ to an affine chart, and embed it in some $\A^d$, so that
it is defined by a set of equations
$$
f_1(u_1,\dots,u_d) = \dots = f_r(u_1,\dots,u_d) = 0.
$$
Then $\J_X$ is locally defined, in this chart, by the
$d-n$ minors of the Jacobian matrix $(\de f_j/\de u_i)$.
Let $S \subset X$ be subscheme defined by $\J_X$. Note that
$S$ is supported precisely over the singular locus of $X$.

We decompose
$$
X_\infty \setminus S_\infty = \bigsqcup_{e=0}^{\infty} X_\infty^e,
\where
X_\infty^e :=  \{ \g\in X_\infty \mid \ord_\g(\J_X) = e \},
$$
and let $X_{m,\infty} := \psi_m(X_\infty)$ and
$X_{m,\infty}^e := \psi_m(X_\infty^e)$,
where $\psi_m \colon X_\infty \to X_m$ is the truncation map.
Also, let
$$
X_\infty^{\le e} : = \{ \g\in X_\infty \mid \ord_\g(\J_X) \le e \}
\and
X_{m,\infty}^{\le e} := \psi_m(X_\infty^{\le e}).
$$
We will need the following geometric lemma on the fibers of the
truncation maps. A weaker version of this property was proven by Denef
and Loeser in \cite[Lemma 4.1]{DL}; the sharper stated here
is taken from \cite[Proposition~4.1]{EM}.

\begin{lem}[\protect{\cite{DL,EM}}]
\label{lemma:fibers-of-pi}
For $m \ge e$, the morphism
$X_{m+1,\infty}^e \to X_{m,\infty}^e$ is a piecewise trivial
fibration with fibers isomorphic to $\A^n$.
\end{lem}

\begin{defi}
\label{quasicylinder} 
A set $C$ in $X_\infty$ that is not contained in $S_\infty$
is a \emph{quasi-cylinder} if there is a cylinder $C'$ 
such that $C \cap (X_\infty\setminus S_\infty)=
\overline{C'}\cap (X_\infty\setminus S_\infty)$.
  \end{defi}


Let $C$ be an irreducible closed quasi-cylinder.
By Lemma~\ref{lemma:fibers-of-pi},
we can define the codimension of $C$ in $X_\infty$.
Indeed, for a cylinder $C'$ such that $C \cap (X_\infty\setminus S_\infty)=
\overline{C'}\cap (X_\infty\setminus S_\infty)$,
one can check that the codimension of
$$
{C'}_m^{\le e} := \psi_m(C') \cap X_{m,\infty}^{\le e}
$$
inside $X_{m,\infty}^{\le e}$ stabilizes for $m\gg e$
(this is done in detail in Section~5 of \cite{EM}), 
and thus we can define
$$
\codim(C',X_\infty) := \codim({C'}_m^{\le e}, X_{m,\infty}^{\le e})
\for m \gg e.
$$
Here, we note that the condition ${C'}_m^{\le e}$ is not empty
is equivalent to the condition
$e \ge \ord_\a(\J_X)$, where $\a\in C'$ is  the generic point,
and therefore ${C'}_m^{\le e}$ is indeed nonempty
for $e \gg 0$, since $C' \not\subseteq S_\infty$. 
Observe that ${C'}_m^{\le e}$ is open, hence dense,
in $\psi_m(C')$.
Note also that 
$\codim(C',X_\infty) = \codim(\overline{C'},X_\infty)$, 
since $\psi_m(C') \subseteq \psi_m(\overline{C'})
\subseteq \overline{\psi_m(C')}$.
Then we define
$$
\codim(C,X_\infty) := \codim(C',X_\infty).
$$
Notice that this definition does not depend on the choice of $C'$. 
Indeed if $C''$ is another cylinder such that
$C \cap (X_\infty\setminus S_\infty)=
\overline{C''}\cap (X_\infty\setminus S_\infty)$, then
we have $\overline{C'} \cap (X_\infty\setminus S_\infty)=
\overline{C''}\cap (X_\infty\setminus S_\infty)$. This implies that 
the symmetric difference between $\overline{C'}$ and $\overline{C''}$
is contained in $S_\infty$. We deduce that
$\codim(\overline{C''},X_\infty) = \codim (\overline{C'},X_\infty)$
and hence $\codim(C'',X_\infty) = \codim (C',X_\infty)$. 

\begin{rmk}
If $s=\codim(\psi_m(C), X_{m,\infty})$ (for $m\gg 0$)
is the codimension of $C$ as defined above,
and $r$ is the maximal length of a sequence
$C=C_0\subsetneq C_1\subsetneq \cdots \subsetneq C_r=X_\infty$
of irreducible closed subsets of $X_\infty$, then
we have the inequality
$$
r\leq s. 
$$
This inequality can be a strict inequality, see for instance \cite[Example~2.8]{IR}.
The inequality can be seen as follows: from the  sequence,
$$C=C_0\subsetneq C_1\subsetneq \cdots \subsetneq C_r$$
of irreducible closed subsets of $X_\infty$,
we have the sequence
$$\overline{\psi_m(C)}=\overline{\psi_m(C_0)}
\subsetneq \overline{\psi_m(C_1)}\subsetneq
\cdots \subsetneq \overline{\psi_m(C_r)}$$
for $m\gg 0$, since $C_i=\varprojlim\overline{\psi_m(C_i)}$.
\end{rmk}

For a non-singular variety \( X \), every component of a cylinder is fat
(\cite{ELM}), but this is no longer true in the singular
case, as the following example show.

\begin{eg}
  Let \( X \) be the Whitney Umbrella, i.e. a hypersurface in
  \( \C^3 \) defined by \( xy^2-z^2=0 \).
  For \( m\geq 1 \),
  let \( \alpha_{m}\colon \C[x,y,z]/(xy^2-z^2)\to \C[t]/(t^{m+1}) \)
  be the \( m \)-jet defined by \( \a_{m}(x)=t, \a_{m}(y)=0, \a_{m}(z)=0 \).
  Then, the cylinder \( \psi_m^{-1}(\a_{m}) \) is contained in \(
  \sing(X)_{\infty} \), where \( \sing(X)= (y=z=0) \).
\end{eg}

\begin{prop}
\label{finitecomponent}
  Let $X$ be a reduced scheme.
  The number of irreducible components of a cylinder on $X_\infty$ is finite.
\end{prop}

\begin{proof}
  First, we show the number of irreducible
components of a cylinder that are not contained in
  $S_\infty$ is finite.
  Let \( C=\psi_{m}^{-1}(\Sigma) \) be a cylinder over a constructible set \(
  \Sigma \subset X_{m} \).
  Let \( \varphi\colon Y\to X \) be a resolution of the singularities of \(
  X \), and assume that $\varphi$ is an isomorphism away from $S$.
    As $\varphi$ is isomorphic away from $S$,
by the valuative criteria for the properness, the generic point of
  each component of \( C \) not contained in $S_\infty$ can be lifted to the
  generic point of a component of  the cylinder
  \(\varphi^{-1}_{\infty}(C)=(\psi_{m}^Y)^{-1}\varphi_{m}^{-1}(\Sigma)\).
  The finiteness of the  components of \( C \) not contained in $S_\infty$
   follows from the finiteness of the components of
   \( (\psi_{m}^Y)^{-1}\varphi_{m}^{-1}(\Sigma) \), as  \(
   Y \) is non-singular.

   Now, to prove the proposition, we use  induction of the dimension.
   If \( \dim X=0 \), then the assertion is trivial,
since $X_\infty\simeq X$  is a finite points set.
   If \( \dim X=n \geq 1 \) and assume that the assertion is true for
   a reduced scheme of dimension \( \leq n-1\).
   Let \( F_{1},\ldots,F_{r} \) be the  irreducible components of
   \( C \) not contained in $S_\infty$.
   Let $$C':=
  \overline{\psi_m\Big(C\setminus
   \bigcup_{i}F_{i}\Big)}$$
   be the closure in   \( \Sigma\cap S_{m} \).
   Then, every irreducible component \( F \) of
   \( C \) other than $F_i$'s  is contained in \( (\psi_m^S)^{-1}(C')\).
      As \( F \) is an irreducible component of \( C=\psi_{m}^{-1}(\Sigma) \)
  which contains  \( (\psi_{m}^{S})^{-1}(C') \),
   \( F  \) is also an irreducible component of the cylinder
   \( (\psi_{m}^{S})^{-1}(C') \) of a lower
   dimensional variety \( S \).
   By  the
   induction hypothesis, we obtain the assertion of our proposition.
   \end{proof}

The second part of the following corollary also appears
as \cite[Lemma~5.1]{EM}.

\begin{cor}\label{thin} Every  fat irreducible
component of a cylinder is a quasi-cylinder, and
every thin component of a cylinder is contained in $S_\infty$.
\end{cor}

\begin{proof}
  Let $C_1,\ldots C_r $ be the irreducible components of a cylinder $C$.
  As $C_i=\varprojlim_m\overline{\psi_m(C_i)}$,
 for $m\gg 0$ it follows that
$${\psi_m(C_i)}\not\subset\overline{\psi_m(C_j)}\ \ \
 \for i\neq j.$$
 Then the non-empty open subset
$C_i\setminus \big( \bigcup_{j\neq i}\psi_m^{-1}(\overline{\psi_m(C_j)})\big)$
of $C_i$ is a cylinder, therefore if $C_i$ is fat,
then  $C_i$ is an irreducible quasi-cylinder.

For the second assertion,   assume $C_1$ is thin.
 If $C_1$ is not contained in $S_\infty$, there is $e>0$ such that
 $X_\infty^{\le e}\cap C_1\neq \emptyset$ which is open in $C_1$.
 Let
 $$U:=X_\infty^{\le e}\setminus
\Big(\bigcup_{j\neq 1}\psi_m^{-1}(\overline{\psi_m(C_j)}) \Big) \and
U_m:=X_{m,\infty}^{\le e}\setminus\Big(\bigcup_{j\neq 1}
 \overline{\psi_m(C_j)}\Big).$$
 Then
 $$U\cap C_1=\psi_m^{-1}(U_m\cap \psi_m(C_1))$$
 is a non-empty open subset of $C_1$ for every $m\gg 0$.
 By Lemma \ref{lemma:fibers-of-pi},
the codimension of $U_m\cap \psi_m(C_1)$  in
 $X_{m, \infty}$ is bounded.
 But it is in contradiction with the fact that $C_1$
is thin by \cite[Lemma 3.7]{must01}.
 \end{proof}

The truncation morphisms induce morphisms
$$
f_m \colon Y_m \lra X_{m,\infty} \subseteq X_m.
$$
Indeed, the inclusion $f_m(Y_m) \subseteq X_{m,\infty}$ is implied by
the fact that $Y_{m,\infty} = Y_m$, which follows by the smoothness of $Y$.
An important ingredient of the  proof for our codimension formula,
as well as the key ingredient
for the change-of-variable formula in motivic integration, is
the following geometric statement on the fibers of these morphisms,
due to Denef and Loeser \protect{\cite[Lemma~3.4]{DL}}.
For a more precise statement of the following lemma,
we refer to \cite[Theorem~6.2 and Lemma~6.3]{EM}.

\begin{lem}[\cite{DL}]
\label{lemma:fibers-of-f_m}
Let $\g \in Y_\infty$ be any arc such that $\t := \ord_\g(\^K_{Y/X}) < \infty$.
Then for any $m \ge 2\t$, letting $\g_m = \psi_m^Y(\g)$, we have
$$
f_m^{-1}\big(f_m(\g_m)\big) \cong \A^\t.
$$
Moreover, for every $\g'_m \in f_m^{-1}(f_m(\g_m))$ we have
$\p^Y_{m,m-\t}(\g_m) = \p^Y_{m,m-\t}(\g'_m)$,
where $\p^Y_{m,m-\t}: Y_m \to Y_{m-\t}$ is the canonical truncation morphism.
\end{lem}

We obtain the following results.

\begin{thm}\label{thm:valuations}
Let $f:Y\to X$ be a resolution of the
singularities such that $E$ appears as a smooth divisor on $Y$.
Then, $W(E,q)$ is a quasi-cylinder of $X_\infty$
of codimension
$$
\codim(W(E,q),X_\infty) = q\.( \ord_E(\^K_{Y/X}) + 1),
$$
\end{thm}

\begin{proof}
Let
$$
\Cont^q(E)^0 := \{\g \in Y_\infty \mid \ord_\g(E) = q, \,
\ord_\g(\ov{\Ex(f) \setminus E}) = 0\}.
$$
This is an open subset of $\Cont^q(E)$. Note that
$\t := \ord_\g(\^K_{Y/X})$ is the same for all $\g \in \Cont^q(E)^0$.
Then, by Lemma~\ref{lemma:fibers-of-f_m}, one can see that
$$
f_m^{-1}\big(f_m\big(\psi^Y_m(\Cont^q(E)^0)\big)\big) = \psi^Y_m(\Cont^q(E)^0)
\fall m \gg 0.
$$
The fact that $\overline{f_\infty(\Cont^q(E)^0)}$
is a quasi-cylinder in $X_\infty$
follows by the equalities $f_m \o \psi^Y_{m} = \psi_{m} \o f_{\infty}$
 and the fact that $\Cont^q(E)^0$ is a cylinder in $Y_\infty$.

 Let $\ord_E\^K_{Y/X} =k$ and $e=q\.k$.
 Then $\psi_m^Y(\cont^q(E))\subset \psi_m^Y(\cont^e(\^K_{Y/X}))$.
 By Lemma~\ref{lemma:fibers-of-f_m}, the morphism
 $\psi_m^Y(\cont^e(\^K_{Y/X}))\to f_m\big(\psi_m^Y(\cont^e(\^K_{Y/X}))\big)$
 induces a morphism
 $$
 \psi_m^Y(\cont^q(E))\to f_m\big(\psi_m^Y(\cont^q(E))\big)
 $$
 with irreducible $e$-dimensional fibers for $m\gg e$.
 Note that $\psi_m^Y(\cont^q(E))$ is an irreducible closed subset of
 codimension $q$ in $Y_m$. Then
$$
\dim \big( f_m\big(\psi_m^Y(\cont^q(E))\big)\big)
=\dim \big(\psi_m^Y(\cont^q(E))\big)-e=(m+1)n-q-e=(m+1)n-q(k+1).
$$
The formula on codimension follows.
\end{proof}

\begin{prop}\label{cyl-con}
  The valuation corresponding to a fat irreducible component of a cylinder
is  the valuation corresponding to an irreducible component of a contact locus.
\end{prop}

\begin{proof}
Let $C \subset X_\infty$ be any fat irreducible component of a cylinder.
 For every $k \in \N$, define
$$
\aa_{k!} := \{f \in \O_X \mid \val_C(f) \ge k! \}
\and
B_{k!} := \{ \g \in X_\infty \mid \ord_\g(\aa_{k!}) \ge k! \}.
$$
Note that we have a chain of inclusions
$$
B_{n!} \supseteq B_{(n+1)!} \supseteq \dots \supseteq C.
$$
For each $k$, let $c_k$ be the smallest
codimension  in $X_\infty$ of irreducible components of $B_{k!}$
containing $C$,
and let $n_k$ be the number of such components.  Since
$$
c_k \le c_{k+1} \le \codim(C,X_\infty) < \infty
\and
$$
$$
 \mbox{if} \  c_k=c_{k+1},\ \mbox{then}\ \   n_k \ge n_{k+1} > 0,
$$
the sequences $\{c_k\}$ and $\{n_k\}$ stabilize.
Therefore we find a closed subset $W \subset X_\infty$ containing $C$
and equal to an irreducible component of $B_{k!}$ for every $k \gg 0$.
We clearly have $\val_C \ge \val_W$ on regular functions
because of the inclusion $C \subseteq W$.
Conversely, let $h \in \O_X$ be an arbitrary nontrivial element.
We can arrange to have $\val_C(h^m) = k!$
for some $m,k \in \N$. This means that $h^m \in \aa_{k!}$,
hence we have $\val_W(h^m) \ge k!$ by the definition of $W$.
Then we conclude that
$$
\val_W(h) = \frac{\val_W(h^m)}m \ge \frac{\val_C(h^m)}m = \val_C(h).
$$
\end{proof}

\section{Determination of a divisorial valuation by finite data}

Let $X = \Spec A$ be an affine variety, and
let \( v \) be a divisorial valuation over $X$, i.e.,
\( v=q \val_{E} \) for \( q\in \N \) and a divisor \( E \)
 over \( X \).
 For a subset \( V\subset X_{\infty} \) we denote the set of
 fat arcs in \(V\) by \( V^o \).

\begin{lem}
\label{blowup}
With the above notation, let
\( x_{1},\ldots, x_{m} \) be elements in \( A \), and denote by
\[ \varphi\colon Y=\Spec A\left[ \frac{x_{2}}{x_{1}},\ldots,
\frac{x_{m}}{x_{1}}\right] \to X = \Spec A\]
the canonical birational morphism.
If \[ \S=\overline{\bigcap_{j=1}^r \left(\cont ^{v_{j}}
(f_{j})\right)^o} \] is an irreducible subset in $Y_\infty$
for some $f_1, \ldots ,f_r \in A\left[ \frac{x_{2}}{x_{1}},\ldots,
\frac{x_{m}}{x_{1}}\right]$ and $v_1,\dots,v_r \in \N$,
then
\begin{equation}\label{eq:3}
 \overline{\varphi_{\infty}(\S)}=
\overline{\left(\bigcap_{j=1}^r \cont ^{v'_{j}}
(f'_{j})\right)\cap \left(\bigcap_{i=1}^m
\cont^{p_{i}}(x_{i})\right)^o}
\end{equation}
for some \( v'_{j}, p_{i} \in \N \) and \( f'_{j}\in A \).
\end{lem}

\begin{proof}
  First we define \( v'_{j}, p_{i} \) and \( f'_{j} \)
  \( (i=1,\ldots, m, j=1,\ldots, r) \).
  Let \( \tilde \alpha \) be the generic point of \( \S \) and
  \( \alpha = \varphi_{\infty}(\tilde \alpha) \).
  Let \( p_{i}=\ord_\a(x_{i})\geq 0 \).
  As \( \alpha \) has the lifting \( \tilde\alpha \) on
  \( Y \), we have \( \ord_ \a(x_{i})-\ord_\a(x_{1})
  = \ord_{\tilde\a}(x_{i}/x_{1})\geq 0 \), which means
  \( p_{1}\leq p_{i} \) for every \( i=2,\ldots,m \).

  Now for \( f_{j}\in A\left[ \frac{x_{2}}{x_{1}},\ldots,
\frac{x_{m}}{x_{1}}\right] \), let the minimal \( a \)
such that \( f_{j}x_{1}^a\in A \) be \( a_{j}  \) and
let \( f'_{j}=f_{j}x_{1}^{a_{j}} \).
  Next let \( v'_{j}=a_{j}p_{1}+v_{j} \).
  Then it is clear that
  \[ \alpha\in \left(\bigcap_{j=1}^r \cont ^{v'_{j}}
(f'_{j})\right)\cap \left(\bigcap_{i=1}^m
\cont^{p_{i}}(x_{i})\right)^o. \]
  Therefore, the inclusion of the left side of \eqref{eq:3}
in the right side follows.
  For the converse inclusion,
  take any arc \( \beta\in \big(\bigcap_{j=1}^r \cont ^{v'_{j}}
(f'_{j})\big)\cap \big(\bigcap_{i=1}^m
\cont^{p_{i}}(x_{i})\big)^o \).
Then, by the condition \( p_{1}\leq p_{i} \) \( (i>1) \), \( \beta   \)
has the
lifting \( \tilde\beta \) on Y.
  Hence, \( \ord_{ \tilde\b}(f_{j})=\ord_\b(f'_{j})
-\ord_ \b(x_{1}^{a_{j}})=v_{j} \) which implies that
\( \tilde\beta\in \Sigma \).
\end{proof}

\begin{thm}
\label{finite}
Let $X = \Spec A$ be an affine variety.
  For any divisorial valuation \( v \) over \( X \),
 there are  finite number of functions \( f_{1},\ldots,f_{r}
  \in A \) and positive integers \( v_{1},\ldots,v_{r} \) such that
  \[ W(v)=\overline{\bigcap_{i=1}^r(\cont ^{v_{i}}(f_{i}))
  ^o}. \]
  Equivalently,
we have $W(v)\sim \bigcap_{i=1}^r\cont ^{v_{i}}(f_{i})$.
\end{thm}

\begin{proof}
  For a divisor \( E \) over \( X \) such that \( v=q \val_{E} \),
  there exists a resolution of singularities
  \( \varphi:Y\to X \) such that \( E \) appears
  as a smooth divisor on \( Y \).
  Then, \( \varphi \) is obtained by blowing up an
  ideal \( I \) in \( A \).
  Let \( g_{1},\ldots,g_{m} \) be  generators of
  \( I \).
  Then \( Y \) is covered by affine open subsets
  \( U_{i}=\Spec A\left[\frac{g_{1}}{g_{i}},\ldots,
  \frac{g_{m}}{g_{i}}\right]
  , (i=1,\ldots,m) \).
  Take an \( i \) such that \( U_{i}\cap E\neq \emptyset \).
  Let this \( i \) be 1.
  Let \( B=A\left[\frac{g_{2}}{g_{1}},\ldots,
  \frac{g_{m}}{g_{1}}\right]\) and consider the restricted
  morphism
  \[ \varphi\colon U_{1}=\Spec B\to X=\Spec A. \]
  Let the defining ideal of \( E \) on \( U_{1} \) be
  \( J=(h_{1},\ldots,h_{s}) \).
  Then the blowing-up \( \psi:U'_{1}\to U_{1} \) of the ideal \( J \) is
  an isomorphism, because it is the blow up with the center
  Cartier divisor \( E \).
  Hence \( U'_{1} \) is non-singular and it
   is covered by charts \( V_{j}=
  \Spec B \left[\frac{h_{1}}{h_{j}},\ldots,
  \frac{h_{s}}{h_{j}}\right]\) \( (j=1,\ldots,s) \).
  Take \( j \) such that \( V_{j}\cap E\neq \emptyset \).
  Let this \( j \) be 1.
  Then, on \( V_{1} \) the divisor \( E \) is a principal divisor
  defined by \( h_{1} \).
  Therefore, on \( V_{1} \), we have
  \( W_{V_{1}}(v)=\overline{(\cont^q(h_{1}))^o} \).
  Applying Lemma \ref{blowup} for
  \[\psi\colon V_{1}=\Spec B \left[\frac{h_{2}}{h_{1}},\ldots,
  \frac{h_{s}}{h_{1}}\right]\to U_{1}=\Spec B,\]
and the property in Remark~\ref{rmk:I3},  we obtain that
  \[ W_{U_{1}}(v)=\overline{\psi_{\infty}(W_{V_{1}}(v))}=
  \overline{(\cont^q(h_{1})\cap \cont^{p_{2}}(h_{2})..
  \cap \cont^{p_{s}}(h_{s}))^o  } \] for some
  \( p_{2},\ldots,p_{s} \).
  Next, applying the lemma and the remark again to
  \[ \varphi: U_{1}=\Spec B\to X=\Spec A, \]
  we obtain
  \[ W_X(v)=\overline{\varphi_{\infty}(W_{U_{1}}(v))  }\]
  \[=\overline{(\cont^{q'}(h'_{1})\cap \cont^{p'_{2}}(h'_{2})\cap\ldots
  \cont^{p'_{s}}(h'_{s})\cap\cont^{q_{1}}(g_{1})\cap\ldots\cap
  \cont^{q_{m}}(g_{m}))^o    } \]
  for some \( h'_{j}\in A \), \( q', p'_{i},q_{l}\in \N \),
(\(1\le j \le s\),
 \( i\ge 2\), and \( 1 \le l \le m) \).
\end{proof}

\begin{rmk}
 By this theorem, it follows that $W(v)$ is the unique fat component of the
quasi-cylinder $\overline{\bigcap_{i=1}^r\cont ^{v_{i}}(f_{i})}$.
As a thin component of a quasi-cylinder is in $S_\infty$
by Corollary \ref{thin}, we see that the inclusion \eqref{eq:inclusion}
is indeed an equality off $S_\infty$, which implies
Theorem \ref{thm:intro:finite}.
  The fact that $W(v)$ is a quasi-cylinder, which was
proven in Theorem \ref{thm:valuations}, also
follows immediately from this theorem.
 \end{rmk}

\begin{defi} Under the situation of Theorem \ref{finite},
we call \( v(f_{1})=q_{1},\ldots, v(f_{r})=q_{r} \)
the {\it generating conditions}, or simply the {\it generators}, of
\( v \).
\end{defi}

 Actually, the conditions in the above definition
``generate'' or ``determine''  the divisorial valuation
  \( v \)  in the  following sense.

\begin{cor}\label{cor:finite-generators}
\label{cfinite}
  Let \( v \) be a divisorial valuation over \( X=\Spec A \).
  Then there exist \( f_{1},\ldots,f_{r} \in A \) and \(
  v_{1},\ldots,v_{r}\in \N \) such that for every \( f\in
  A\setminus \{0\} \)
  \[ v(f)=\min\{v'(f) \mid \text{$v'$ is a divisorial valuation over
 $X$ such that $v'(f_{i})=v_{i}$} \}. \]
\end{cor}

\begin{proof}
  Let \( v(f_{i})=v_{i}\)  \(
 (i=1,\ldots,{r}) \)
  be generators of \( v \).
  For every $f \in A \setminus \{0\}$, let \( v_{0}(f) \) be
  \[ \min\{v'(f) \mid \text{$v'$ is a divisorial valuation over
 $X$ such that $v'(f_{i})=v_{i}$}\}. \]
  Then it is clear that \( v_{0}(f)\leq v(f) \) for every \( f\in
  A\setminus \{0\} \) by the definition of \(
  v_{0} \).
  On the other hand for an arbitrary \( f\in A\setminus\{0\} \),
  take a divisorial valuation \( v' \)  over \( X \)
(depending on \( f \)) such that
  \( v'(f)=v_{0}(f) \) and \( v'(f_{i})=v_{i} \) \( (i=1,\ldots,r) \).
  Let \( \beta \) be the generic point of \( W(v') \).
  Then \( \beta\in \overline{\bigcap_{i=1}^r(\cont ^{v_{i}}(f_{i}))
  ^o}=W(v)  \).
  The inclusion \( W(v')\subset W(v) \) yields the inequality
  \( v'\mid _{A}\geq v\mid _{A} \).
    Therefore, in particular, we have \( v_{0}(f)=v'(f)\geq v(f) \).
\end{proof}

\begin{rmk}\label{rmk:MacLane}
Suppose that $K$ is the function field of a $n$-dimensional
variety and $x$ is an element purely transcendental over $K$. Suppose
that $v$ is a divisorial valuation for $K(x)$ with residue field L.
Let  $w$ be the restriction of $v$ to $K$. Denote by $L'$ the residue field
of $w$. MacLane's theorem says that one can construct $v$ from $w$
by knowing the values of $v$ on a sequence of key polynomials in $K[x]$
(see \cite{ML1}).
If the sequence is infinite, then $L$ is algebraic over $L'$. Then the
transcendence degree of $L$ over $\C$ would be at most $n-1$ and $v$  can
not be  divisorial. Hence $v$ can be determined from $w$
and the values of $v$ on a finite number of key polynomials. Inductively,
we see that we can determine a divisorial valuation $v$ on a finitely
generated purely
transcendental field over $\C$ by knowing the values of $v$
on  just finitely many
elements in the function fields.
\end{rmk}

\begin{eg}[Divisorial valuations over \( \C^2 \)]
  We give some explicit constructions of divisorial valuations
of $\C(x,y)$ in terms of blow-ups of \( \C^2 \).
We will denote by \( v_{1},\dots,v_{4},v_{5} \) and \( v'_{5} \)
the valuations \(
  \val_{E_{1}},\dots,\val_{E_{4}}\), \(\val_{E_{5}} \) and \( \val_{E'_{5}} \)
 associated to the following exceptional divisors
\( E_{1},\dots,E_{4},E_{5} \) and   \( E'_{5} \).
\begin{enumerate}
 \item[(i)]
   Let \( E_{1} \) be the exceptional divisor of the blowing-up
   \( \varphi_{1}:X_{1}\to \C^2 \) with the center \( (x, y)=(0,0) \).
   Then, the generators of \(v_{1}\)
are \( v_{1}(x)=v_{1}(y)=1 \) and \( v_{1} \)
is a toric valuation.
 \item[(ii)]
   Let \( E_{2} \) be the exceptional divisor of the blowing-up
   \( \varphi_{2}:X_{2}\to X_{1} \) with the center \(( x_{1}, y_{1})=
   (0,0) \), where \( x_{1}=x, y_{1}=y/x \).
   Then, the generators of \(v_{2}\) are \( v_{2}(x)=1 \),
\( v_{2}(y)=2 \) and
\( v_{2} \) is a toric
   valuation.
 \item[(iii)]
   Let \( E_{3} \) be the exceptional divisor of the blowing-up
   \( \varphi_{3}:X_{3}\to X_{2} \)  with the center \( (x_{2},y_{2})=
   (0,0) \), where \( x_{2}=x_{1}/y_{1} \), \( y_{2}=y_{1} \).
   Then, the generators of \(v_{3}\) are   \( v_{3}(x)=2 \), \( v_{3}(y)=3 \)
    and \( v_{3} \) is a toric
   valuation.
 \item[(iv)]
   Let \( E_{4} \) be the exceptional divisor of the blowing-up
   \( \varphi_{4}:X_{4}\to X_{3} \)  with the center \( (x_{3},y_{3})=
   (0,\lambda) \) \( (\lambda\in \C\setminus \{0\}) \),
   where \( x_{3}=x_{2} \), \( y_{3}=y_{2}/x_{2} \).
   Then, the generators of \( v_{4} \) are \( v_{4}(x)=2 \), \( v_{4}(y)=3 \)
    and
   \( v_{4}(y^2-\lambda x^3)=7 \).
 \item[(v)]
    Let \( E_{5} \) be the exceptional divisor of the blowing-up
   \( \varphi_{5}:X_{5}\to X_{4} \)  with the center \( (x_{4},y_{4})=
   (0,0)\), where \( x_{4}=x_{3} \), \( y_{4}=(y_{3}-\lambda)/x_{3} \).
   Then, the generators of \( v_{5} \) are \( v_{5}(x)=2 \), \( v_{5}(y)=3 \)
   and
   \( v_{5}(y^2-\lambda x^3)=8 \).
 \item[(v')]
     Let \( E'_{5} \) be the exceptional divisor of the blowing-up
   \( \varphi'_{5}:X'_{5}\to X_{4} \)  with the center \( (x_{4},y_{4})=
  (0, \mu)\) \( (\mu\in \C\setminus \{0\}) \).
    Then, the generators of \( v'_{5} \) are \( v'_{5}(x)=2 \), \(
    v'_{5}(y)=3 \),
    \( v'_{5}( y^2-\lambda x^3)=7 \) and \( v'_{5}( y^3-\lambda x^3y-\mu x^5)
    = 11 \).
\end{enumerate}
\end{eg}

\begin{eg}[Divisorial valuations over \( \C^3 \)]
In \cite{cgp} V. Cossart, C. Galindo and O. Piltant showed an example of divisorial valuation $v$ with the center at a
 non singular point $x\in X$ of a three dimensional variety such that the graded ring $\Gr_v(\O_{X,x})$ associated to this valuation is not Noetherian.
But this valuation is also ``finitely generated" in our sense.

In \cite{cgp} the valuation $v$ is obtained as follows:
Let $R$ be a regular local ring of dimension three and let $X_0=\Spec R$.
Let $X_1\to X_0$ be the blow up at the closed point of $X_0$, $E_1=\Proj \C[x,y,z] $ the exceptional divisor and $C$ the cubic in $E_1$ defined by the equation $x^2z+xy^2+y^3=0$.
Take the blow up $X_2\to X_1$ at the smooth point of $C$ with the homogeneous coordinates $(1:0:0)$.
Then, for $n\geq 10$ construct a sequence  $X_n\to X_{n-1}\to\cdots \to X_2$, where each $X_{i+1}\to X_{i}$ is the blow up at the intersection of the strict transform of $C$ and the  exceptional divisor of $X_{i}\to X_{i-1}$.
For simplicity, we assume that $R$ is the local ring at the origin of $X=\Spec\C[x,y,z]$.
As $\Gr_v(R)$ is not Noetherian, the graded ring $\Gr_v(\C[x,y,z])$ is also non-Noetherian.
But on $X$ we have finite "generators": $v(x)=n-1$, $v(y)=n$, $v(z)=n+1$ and $v(x^2z+xy^2+y^3)=4(n-1)$.
\end{eg}

\section{Characterization of toric valuations}

Using the result of \cite{ELM}, Theorem~\ref{thm-intro:monomial-smooth}
follows immediately. Although it is a particular case of the
general statement for toric varieties, it is instructive to
give the proof independently.

\begin{proof}[Proof of Theorem~\ref{thm-intro:monomial-smooth}]
Let $X = \C^n$, and consider the set
$$
C := \{ \g \in X_\infty \mid \ord_\g(x_i) \ge a_i \}.
$$
It is immediate to check that
this is an irreducible closed cylinder of codimension $\sum a_i$
in $X_\infty$, and in fact it coincides with the maximal divisorial set
associated with the toric valuation given by assigning weights
$a_i$ to the coordinates $x_i$.
Thus the statement follows as long as we show that
 $W(v) \subset X_\infty$
is equal to $C$. By the first condition in \eqref{eq-intro:monomial-smooth},
we have $W (v)\subseteq C$, whereas the second condition
and Theorem~\ref{thm:valuations} imply that
$\codim(W,X_\infty) \le \codim(C,X_\infty)$. We conclude that $W (v)= C$
by the irreducibility of $C$.
\end{proof}

Using notation as in \cite{Ful}, let $M$
be a free abelian group of rank $n \ge 1$,
and let $N = M^*$ be its dual. Fix a maximal
dimensional rational polyhedral cone
$\s \subset N \otimes \R$, and let $\D$ be the
fan consisting of all faces of $\s$.
Denote $R_\s = \s^\vee \cap M$, and let $T(\D) = \Spec\C[R_\s]$
be the corresponding toric variety.
Recall that a toric valuation on a toric variety $X=T(\D)$ is determined
by a nonzero element in $\s\cap N$, and conversely.

\begin{thm}\label{thm-intro:monomial-singular}
Let $v$ be a divisorial valuation over a complex affine toric
variety $X = \Spec\C[R_\s]$ centered at the origin of $X$.
Assume that there is an element $a \in \s\cap N\setminus \{0\}$ such that
\begin{equation}\label{eq-intro:monomial-singular-1}
  \^k_v(X) + q(v) \le \min_{\{x_1,\dots,x_n\}}
\left\{\sum \langle x_i,a \rangle \right\},
  \end{equation}
where the minimum is taken
over linearly independent subsets $\{x_1,\dots,x_n\} \subset R_\s$
of cardinality $n = \dim X$, and that
\begin{equation}\label{eq-intro:monomial-singular-2}
  v(u_j) \ge \langle u_j,a \rangle \fall j
  \end{equation}
for some set of monomial generators $\{u_1,\dots,u_r\}$  of $R_\s$.
Then $v$ coincides with the toric valuation $\val_a$ determined by $a$,
$q(v) = \max\{k \in \N \mid a \in kN \}$, and equality
holds in \eqref{eq-intro:monomial-singular-1} and
\eqref{eq-intro:monomial-singular-2}
\end{thm}

The proof for the singular toric case is similar, but we need to
fix a couple of properties first.
From now on, suppose that $X = T(\D)$ (notation as above).
Fix a nonzero element $a \in \s\cap N$,
and let $\val_a$ be the associated toric valuation.

\begin{lem}\label{lem:Mather-discr-toric-case}
With the above notation, we have
$$
\^k_{\val_a}(X) =
\min_{\{x_1,\dots,x_n\}} \left(\sum \langle x_i,a \rangle - q(a) \right),
$$
where $q(a) = \max\{k \in \N \mid a \in kN \}$, and
the minimum is taken over the set of linearly independent subsets
$\{x_1,\dots,x_n\}$ of cardinality $n$ in $R_\s$.
\end{lem}

\begin{proof}
By replacing $a$ by its primitive element $a/q(a)$,
we can assume that $q(a) = 1$.
We observe that $\wedge^n \Om_X$ is generated by the
torus-invariant forms
$$
{dx_1} \wedge \dots \wedge{dx_n},
$$
as $\{x_1,\dots,x_n\}$ ranges among linearly independent
subsets of cardinality $n$ in $\s^\vee \cap M$.

Take a toric resolution of singularities $f \colon Y \to X$
factoring both through the toric blowup $g \colon \Bl_a(X) \to X$
determined by $a$, and through the Nash blowup $\n \colon \^X \to X$.
Let $E \subset Y$ be the prime divisor that is the proper transform of the
exceptional divisor of $g$, and fix toric invariant coordinates
$(y_1,\dots,y_n)$ in a toric affine open set of $Y$ intersecting $E$.
Note that $E$ is defined by the vanishing of one of these $y_j$.

Note that $\wedge^n\Om_Y \cong \O_Y(K_Y)$ is generated (in this open set)
by the form
$$
{dy_1} \wedge \dots \wedge {dy_n}.
$$
On the other hand, the image of the map
$df \colon f^*\wedge^n\Om_X \to \wedge^n\Om_Y$
is determined, in this open set, by
$$
dx_1 \wedge\dots\wedge dx_n = (x_1\dots x_n)\.
\frac{dy_1}{y_1} \wedge \dots \wedge \frac{dy_n}{y_n}
\quad\text{up to unit},
$$
as $\{x_1,\dots,x_n\}$ ranges as above. The statement follows by the
definition of $\^K_{Y/X}$ as the divisor determined by the image of $df$.
\end{proof}

The next step is to understand the maximal divisorial set  associated to
the toric valuation $\val_a$.
We fix nonconstant monomial generators $u_1,\dots,u_r$ of $R_\s$, and let
$a_j := \langle u_j,a \rangle$.

\begin{lem}\label{lem:W(val_a)}
The maximal  divisorial set  associated with $\val_a$ is equal to
$$
W(\val_a) = \{ \g \in X_\infty \mid \ord_\g(u_j) \ge a_j \}.
$$
\end{lem}

\begin{proof}
The inclusion $\subset$ is trivial.
We prove the opposite inclusion $\supset$.
Let $T \cong (\C^*)^n$ be the torus acting on $X$.
For every $m \ge 1$ we have a commutative diagram
$$
\xymatrix@C=20pt@R=20pt{
T_\infty \ar[d] \ar@{^{(}->}[r] & X_\infty \ar[d] \\
T_m \ar[d] \ar@{^{(}->}[r] & X_m \ar[d] \\
T \ar@{^{(}->}[r] & X,
}
$$
where $T_m \cong \C^{nm} \times (\C^*)^n$ acts on $X_m$
(\cite[Proposition 2.6]{Ish}), and $T_\infty$
is the pro-scheme obtained as the inverse limit of the projective system
$\{T_m \to T_{m-1}\}$.

For every face $\t \preceq \s$ we define
$$
X_\infty(\t) := \{\g \in X_\infty \mid \g(\e) \in \orb(\t) \},
$$
in particular for $\t=\{0\}$
$$
X_\infty(0) := \{\g \in X_\infty \mid \g(\e) \in T \},
$$
where $\e$ denotes the generic point of $\Spec\C(\g)[[t]]$ and
$\C(\g)$ is the residue field of $\g\in X_\infty$.
We obtain a stratification
$$
X_\infty = \bigsqcup_{\t \preceq \s} X_\infty(\t).
$$

Observe that every arc $\g \colon \C[R_\s] \to \C(\g)[[t]]$
determines a semigroup
map
$$
\ord_\g \colon R_\s \to \Z_{\ge 0} \cup \{\infty\}.
$$
If $\g \in X_\infty(0)$, then we actually have
$\ord_\g \colon R_\s \to \Z_{\ge 0}$, and this extends to
a linear map $M \to \Z$. So $\g$ naturally determines an
element $a_\g \in N$, and since $\ord_\g$ is nonnegative
on $R_\s$, we have $a_\g \in \s$.

Now we stratify $X_\infty(0)$ accordingly to the position of $a_\g$
(we treat the other strata $X(\t)$, for $0 \ne \t \preceq \s$,
by induction on dimension) as follows.
For every $a\in \sigma\cap N\setminus \{0\}$ we define
$$T_\infty(a)=\{\g\in X_\infty(0)\mid a_\g=a\}.
$$
Then $T_\infty(a)$ is a $T_\infty$-orbit and
$\overline{T_\infty(a)}=W(\val_a)$.

By \cite[Proposition~4.8]{Ish}, for every two
$a,a' \in (\s \cap N) \setminus \{0\}$ we also have
$$
\ov{T_\infty(a')}\subseteq \ov{T_\infty(a)}
\quad\Longleftrightarrow\quad
a' - a \in \s \cap N
\quad\Longleftrightarrow\quad
\val_{a'}|_{R_\s} \ge \val_{a}|_{R_\s}.
$$

Now take an arbitrary $\g  \in  X_\infty(0)\setminus W(\val_a) $.
Then $\g \in T_\infty(b)$ for some $b \in (\s \cap N) \setminus \{0\}$ for
which $b - a \not\in \s \cap N$, hence we have
$\ord_\g(u_j) < a_j$ for some $u_j$. The statement follows
by iterating the argument on the different pieces $X(\t)$ of the
stratification of $X$ by using \cite[Theorem 4.15]{Ish}.
\end{proof}

\begin{proof}[Proof of Theorem~\ref{thm-intro:monomial-singular}]
The proof goes along the same lines of the one
of Theorem~\ref{thm-intro:monomial-smooth}.
By Lemma~\ref{lem:W(val_a)} and \eqref{eq-intro:monomial-singular-2},
we have $W(v) \subseteq W(\val_a)$. The conclusion follows
by comparing codimension, which are determined respectively by
\eqref{eq-intro:monomial-singular-1}
and Lemma~\ref{lem:Mather-discr-toric-case}.
\end{proof}

\providecommand{\bysame}{\leavevmode \hbox \o3em
{\hrulefill}\thinspace}

\end{document}